\documentclass[11pt]{article}

\usepackage[pdftex]{graphicx}
\usepackage{amssymb}
\usepackage{amsthm}

\makeatletter
\renewcommand\section{\@startsection{section}{1}{\z@}%
                                  {-3.5ex \@plus -1ex \@minus -.2ex}%
                                  {2.3ex \@plus.2ex}%
                                  {\normalfont\large\bfseries}}
\makeatother

\DeclareGraphicsExtensions{.jpg}
\begin{document}

\title{Hamiltonian cubic bipartite graphs}

\author{Misa Nakanishi \thanks{E-mail address : nakanishi@mx-keio.net}}
\date{}
\maketitle

\begin{abstract}
We provide a polynomial time algorithm to determine a cubic bipartite graph has a hamilton cycle or not. \\
keywords : hamilton cycle, cubic bipartite graph, NP-complete
\end{abstract}

\newtheorem{thm}{Theorem}[section]
\newtheorem{lem}{Lemma}[section]
\newtheorem{prop}{Proposition}[section]
\newtheorem{cor}{Corollary}[section]
\newtheorem{rem}{Remark}[section]
\newtheorem{conj}{Conjecture}[section]
\newtheorem{claim}{Claim}[section]
\newtheorem{fact}{Fact}[section]
\newtheorem{obs}{Observation}[section]

\newtheorem{defn}{Definition}[section]
\newtheorem{propa}{Proposition}
\renewcommand{\thepropa}{\Alph{propa}}
\newtheorem{conja}[propa]{Conjecture}

\section{Introduction}
This article takes the hamilton cycle problem and aims to determine by using an algorithm whether a cubic bipartite graph has a hamilton cycle. A cubic bipartite graph is a kind of graph, with no odd cycles and all vertices degree 3. A hamilton cycle is a cycle that visits each vertex exactly once.

\section{Contribution}
The hamilton cycle problem has been the subject of much research for a long time. This article deals with cubic bipartite graphs and resolves this problem by using a polynomial time algorithm.

In 1980, the problem of determining whether a 2-connected cubic bipartite planar graph or a 3-connected cubic bipartite graph has a hamilton cycle was shown to be NP-complete \cite{Akiyama}. Therefore, the problem of determining whether a cubic bipartite graph has a hamilton cycle is NP-complete.

Also, in 1982, it was shown that the problem of determining whether a general grid graph has a hamilton cycle is NP-complete \cite{Itai}. Recently, in 2017, it was shown that the problem of determining whether the line graph of a planar cubic bipartite graph has a hamilton cycle is NP-complete \cite{Munaro}. From these two studies, that a graph is cubic and bipartite is not a necessary condition for the hamilton cycle problem to be NP-complete.

In class NP, the class of problems that can be solved by a polynomial time algorithm is called P, and it is generally said that NP-complete problems are not in P. In this article, we present that one of NP-complete problems is in P and lead that all NP-complete problems are contained in P (see Figure 1).

\includegraphics[width=10cm, trim=0cm 3cm 3cm 2cm, clip=true]{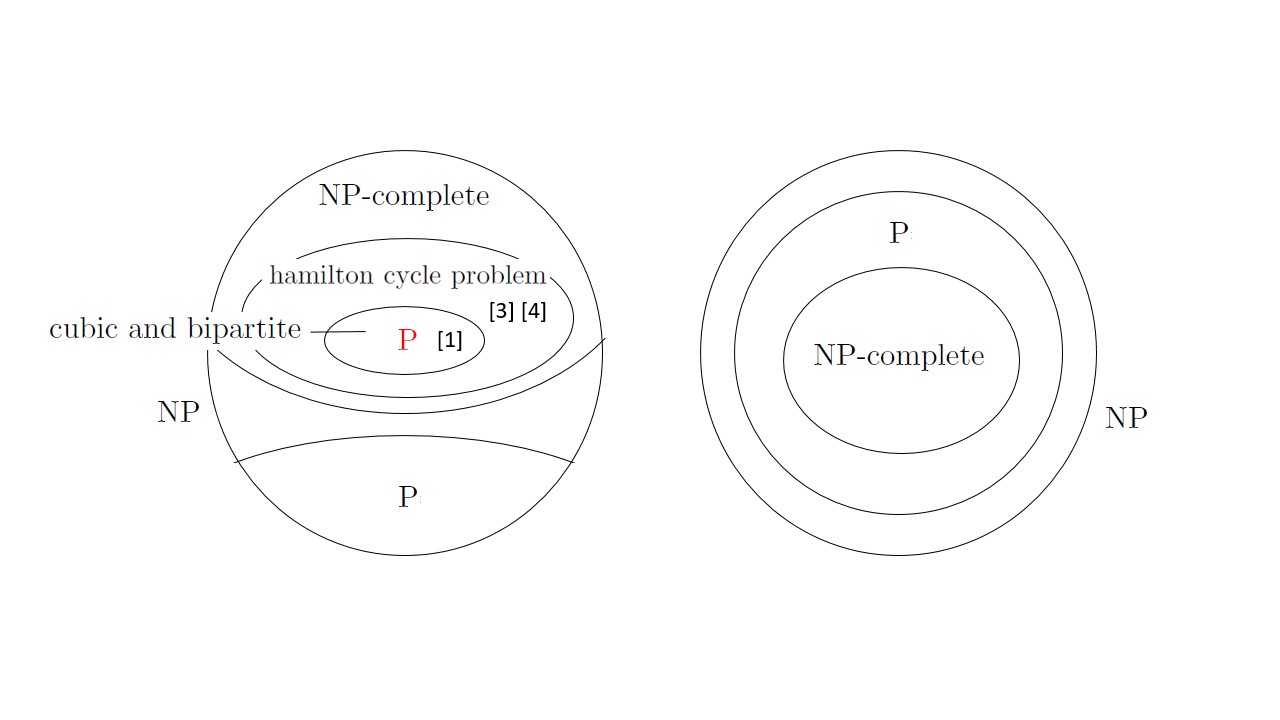}\\
Figure 1: Euler diagram for class NP. In contrast to the conventional diagram on the left (black line), this article shows that the hamilton cycle problem for cubic bipartite graphs is in P (red string), and leads that all NP-complete problems are contained in P as in the diagram on the right. \\

Section \ref{N} introduces the notation and property used in this article. In Section \ref{H}, we show a polynomial time algorithm for solving the hamilton cycle problem for cubic bipartite graphs and prove that the polynomial time algorithm is correct.

\section{Notation and property}\label{N}
We follow the notations presented in \cite{Diestel}. Let $G := (V, E)$ be a graph with the vertex set $V$ and edge set $E$. In this article, a graph is finite, undirected, and simple. For a vertex $v \in V(G)$, the open neighborhood, denoted by $N_G(v)$, is $\{ u \in V(G) \colon\ uv \in E(G) \}$. For an edge $e \in E(G)$, the set of two ends of $e$ is denoted by $V(e)$.   

\begin{propa}[\cite{Petersen}]
Every bridgeless cubic graph has a perfect matching. (Note that the complement of a perfect matching is always a 2-factor.)
\end{propa}

\section{Hamiltonian cubic bipartite graphs}\label{H}

We define the algorithms ${\bf HC}$ and ${\bf DC}$ as follows. \\

\noindent ${\bf HC}$ : Input $G$ as a cubic bipartite graph. \\
(1) Let $M$ be a perfect matching of $G$, and let $F := G - M$. \\ 
(2) Let $M' \subseteq M$ be a maximal subset of $M$ such that for every $m \in M'$, both ends of $m$ are not part of the same cycle of $F$. Let $M'' := M \setminus M'$. \\
(3) Let $K := {\bf DC}(F, M', M'', G)$. Let $M(K) := M \cap E(K)$, $F(K) := E(K) \setminus M(K)$, and $F := F - F(K) + M(K)$. If $F$ is a hamilton cycle of $G$ then go to (6). \\
(4) For $M := E(G) \setminus E(F)$, redefine $M'$ and $M''$ as (2). \\
(5) Repeat (3) - (4) at most $|V(G)|^2$ times. \\
(6) Return $F$. \\

\noindent ${\bf DC}$ : Input $F, M', M''$, and $G$. \\
(1) Let $J := G$. \\
(2) Choose $m \in M' \cup M''$ arbitrarily and let $v_1, v_2 \in V(m)$. (The term mark appears after (4).) \\
\ (*1) If $d_J(v_1) = 3$ and $v_1$ is an end of some marked edge, then for $w \in N_J(v_1) \setminus \{v_2\}$ where $v_1w$ is not marked, let $J := J - v_1w$. \\
\ (*2) If $d_J(v_2) = 3$ and $v_2$ is an end of some marked edge, then for $w \in N_J(v_2) \setminus \{v_1\}$ where $v_2w$ is not marked, let $J := J - v_2w$. \\
(3) Repeat (2) until in (2) all $m \in M' \cup M''$ are chosen. \\
(4) Choose $m \in M' \cup M''$ arbitrarily and let $v_1, v_2 \in V(m)$. \\
\ (*1) If $d_J(v_1) = 2$, then for $w \in N_J(v_1) \setminus \{v_2\}$, mark $v_1w$. \\
\ (*2) If $d_J(v_2) = 2$, then for $w \in N_J(v_2) \setminus \{v_1\}$, mark $v_2w$. \\
(5) Repeat (2) - (4) until in (4) all $m \in M' \cup M''$ are chosen. \\
(6) Choose $m \in M'$ arbitrarily and let $v_1, v_2 \in V(m)$. \\
\ (*1) If $d_J(v_1) = 3$ and $v_2$ is an end of some marked edge, then let $w_1, w_2 \in N_J(v_1) \setminus \{v_2\}$, $x_1 \in N_F(w_1) \setminus \{v_1\}$, and $x_2 \in N_F(w_2) \setminus \{v_1\}$. \\
\ \ (a) If $w_1$ is an end of some $m' \in M'$, and a cycle $C$ such that $v_1w_1 \in E(C)$ exists in $J - v_1w_2Fx_1w_1$, then mark $v_1w_1$ and let $J := J - v_1w_2$. \\
\ \ (b) If $w_2$ is an end of some $m' \in M'$, and a cycle $C$ such that $v_1w_2 \in E(C)$ exists in $J - v_1w_1Fx_2w_2$, then mark $v_1w_2$ and let $J := J - v_1w_1$. \\
\ (*2) If $d_J(v_2) = 3$ and $v_1$ is an end of some marked edge, then let $w_1, w_2 \in N_J(v_2) \setminus \{v_1\}$, $x_1 \in N_F(w_1) \setminus \{v_2\}$, and $x_2 \in N_F(w_2) \setminus \{v_2\}$. \\
\ \ (a) If $w_1$ is an end of some $m' \in M'$, and a cycle $C$ such that $v_2w_1 \in E(C)$ exists in $J - v_2w_2Fx_1w_1$, then mark $v_2w_1$ and let $J := J - v_2w_2$. \\
\ \ (b) If $w_2$ is an end of some $m' \in M'$, and a cycle $C$ such that $v_2w_2 \in E(C)$ exists in $J - v_2w_1Fx_2w_2$, then mark $v_2w_2$ and let $J := J - v_2w_1$. \\
(7) Repeat (2) - (6) until in (6) all $m \in M'$ are chosen. \\
(8) Choose $m \in M'$ arbitrarily and let $v_1, v_2 \in V(m)$. \\
\ (*1) If $d_J(v_1) = 3$, then let $w_1, w_2 \in N_J(v_1) \setminus \{v_2\}$. \\
\ \ (a) If $w_1$ is an end of some $m' \in M'$, and $w_2$ is an end of some $m'' \in M''$, then mark $v_1w_1$, let $J := J - v_1w_2$, and go to (9). \\
\ \ (b) If $w_2$ is an end of some $m' \in M'$, and $w_1$ is an end of some $m'' \in M''$, then mark $v_1w_2$, let $J := J - v_1w_1$, and go to (9). \\
\ (*2) If $d_J(v_2) = 3$, then let $w_1, w_2 \in N_J(v_2) \setminus \{v_1\}$. \\
\ \ (a) If $w_1$ is an end of some $m' \in M'$, and $w_2$ is an end of some $m'' \in M''$, then mark $v_2w_1$, let $J := J - v_2w_2$, and go to (9). \\
\ \ (b) If $w_2$ is an end of some $m' \in M'$, and $w_1$ is an end of some $m'' \in M''$, then mark $v_2w_2$, let $J := J - v_2w_1$, and go to (9). \\
(9) Repeat (2) - (8) until in (8) all $m \in M'$ are chosen. \\
(10) Choose $m \in M'$ arbitrarily and let $v_1, v_2 \in V(m)$. \\
\ (*1) If $d_J(v_1) = 3$, then let $w_1, w_2 \in N_J(v_1) \setminus \{v_2\}$. \\
\ \ (a) Mark $v_1w_1$, let $J := J - v_1w_2$, and go to (11). \\
\ \ (b) Mark $v_1w_2$, let $J := J - v_1w_1$, and go to (11). \\
\ (*2) If $d_J(v_2) = 3$, then let $w_1, w_2 \in N_J(v_2) \setminus \{v_1\}$. \\
\ \ (a) Mark $v_2w_1$, let $J := J - v_2w_2$, and go to (11). \\
\ \ (b) Mark $v_2w_2$, let $J := J - v_2w_1$, and go to (11). \\
(11) Repeat (2) - (10) until in (10) all $m \in M'$ are chosen. \\
(12) Let $K' := J$. \\
(13) Choose $m \in M' \cap K'$ arbitrarily and let $v_1, v_2 \in V(m)$. If $N_{K'}(v_1) \setminus \{v_2\}$ or $N_{K'}(v_2) \setminus \{v_1\}$ is an end of some edge in $M'$, then check $m$. \\
(14) Repeat (13) until in (13) all $m \in M' \cap K'$ are chosen. \\
(15) Choose $m \in M' \cap K'$ arbitrarily and let $v_1, v_2 \in V(m)$. If $N_{K'}(v_1) \setminus \{v_2\}$ or $N_{K'}(v_2) \setminus \{v_1\}$ is an end of some edge in $M''$, then highlight $m$. \\
(16) Repeat (15) until in (15) all $m \in M' \cap K'$ are chosen. \\
(17) Choose a cycle $C \subseteq K'$ arbitrarily. Let $k$ be the number of components in $F \cup (K' - C)$. If $E(C)$ dose not contain any checked edges and $k = 1$, then let $K' := K' - C$. \\
(18) Repeat (17) until in (17) all cycles $C \subseteq K'$ are chosen. \\
(19) Choose a cycle $C \subseteq K'$ arbitrarily. Let $k$ be the number of components in $F \cup (K' - C)$. If $E(C)$ contains some highlighted edges and $k = 1$, then let $K' := K' - C$. \\
(20) Repeat (19) until in (19) all cycles $C \subseteq K'$ are chosen. \\
(21) Choose a cycle $C \subseteq K'$ arbitrarily. Let $k$ be the number of components in $F \cup (K' - C)$. If $k = 1$, then let $K' := K' - C$. \\
(22) Repeat (21) until in (21) all cycles $C \subseteq K'$ are chosen. \\
(23) Let $K'' := K'$. Return $K''$. \\

Let $G$ be a cubic bipartite graph. ${\bf HC}(G)$ is a 2-factor constructed by applying ${\bf HC}$ to $G$. Note that ${\bf HC}(G)$ is not unique and taken arbitrarily.

\includegraphics[width=10cm, trim=6cm 6cm 10cm 3cm, clip=true]{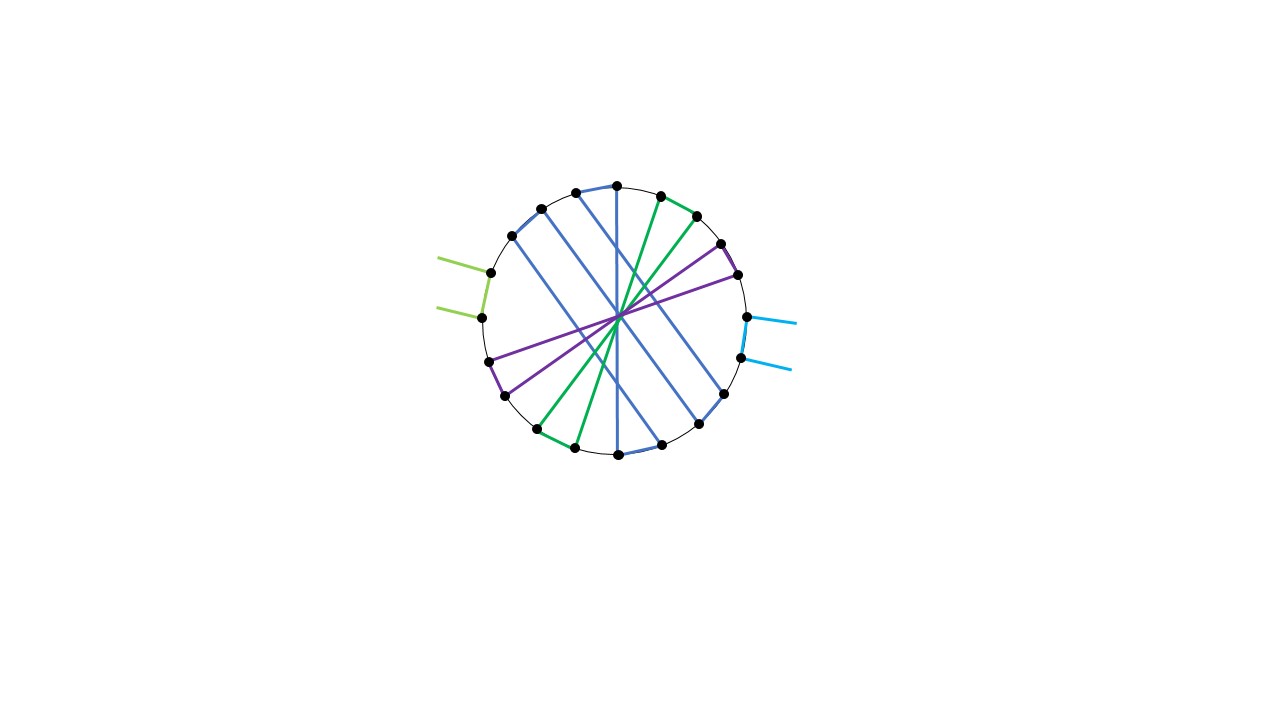}
\begin{center}
Figure 2: Detecting cycles in ${\bf DC}$
\end{center}

\begin{fact}\label{F00}
$K'$ in ${\bf DC}$ (12) is the union of disjoint cycles with alternating edges of $F$ and edges of $M'$ or $M''$, and is the spanning subgraph of $G$.
\end{fact}

\begin{proof}
In ${\bf HC}$, $G$ is a cubic bipartite graph and $F$ is a 2-factor of $G$. In ${\bf DC}$ (2) - (11), if one edge is marked in the cycle of $F$, every other edge is marked from that edge in the cycle of $F$, and unmarked edges are deleted from $J$ ($:= G$). As a result, $K'$ in ${\bf DC}$ (12) is the union of disjoin cycles with alternating edges of $F$ and edges of $M'$ or $M''$, and is the spanning subgraph of $G$ (see Figure 2).
\end{proof}

\begin{fact}\label{23}
$K''$ in ${\bf DC}$ (23) is a minimal subgraph consisting of cycles of $K'$ in ${\bf DC}$ (12) so that $F \cup K''$ has 1 component.
\end{fact}

\begin{proof}
By Fact \ref{F00}, $K'$ in ${\bf DC}$ (12) is the union of disjoint cycles of $G$ and $F \cup K'$ has 1 component. By ${\bf DC}$ (13) - (22), $K''$ in ${\bf DC}$ (23) is a minimal subgraph consisting of cycles of $K'$ in ${\bf DC}$ (12) so that $F \cup K''$ has 1 component.
\end{proof}

\begin{lem}\label{F1}
If $G$ is hamiltonian, then ${\bf HC}(G)$ returns a hamilton cycle.
\end{lem}

\begin{proof}
Apply ${\bf HC}$ to $G$. In ${\bf HC}$ (1) and (2), let $M$ be a perfect matching of $G$, $F := G - M$, $M'$ be the set of edges of $M$ whose ends are not in the same cycle of $F$, and $M'' := M \setminus M'$. In ${\bf HC}$ (3), apply ${\bf DC}$ to $F$, $M'$, $M''$, and $G$. If $F$ is a hamilton cycle, then in ${\bf DC}$ (23), $K'' = \emptyset$. Suppose that $F$ is not a hamilton cycle. By Fact \ref{F00}, $K'$ in ${\bf DC}$ (12) is the union of disjoint cycles with alternating edges of $F$ and edges of $M'$ or $M''$, and is the spanning subgraph of $G$. Indeed, by ${\bf DC}$ (6) and (8), $K'$ in ${\bf DC}$ (12) is constructed by taking a cycle with alternating edges of $F$ and edges of $M'$ with the first priority, and a cycle containing a path with an edge of $M'$, an edge of $F$, an edge of $M'$ in order with the second priority (in particular, ${\bf DC}$ (6) marks the edge of $F$ to close the cycle). (See Figure 3.) By Fact \ref{23}, $K''$ in ${\bf DC}$ (23) is a minimal subgraph consisting of cycles of $K'$ in ${\bf DC}$ (12) so that $F \cup K''$ has 1 component. Indeed, by ${\bf DC}$ (17), (19), and (21), $K''$ in ${\bf DC}$ (23) is constructed by remaining a cycle with alternating edges of $F$ and edges of $M'$ with the first priority, and a cycle containing a path with an edge of $M'$, an edge of $F$, an edge of $M'$ in order with the second priority (in particular, ${\bf DC}$ (17) removes every cycle with alternating edges of $F$ and edges of $M''$). (See Figure 4 and 5.) In ${\bf HC}$ (3), put ${\bf DC}(F, M', M'', G)$ as $K$, and let 
\[M(K) := M \cap E(K),\ F(K) := E(K) \setminus M(K),\]
and 
\begin{equation}
F := F - F(K) + M(K).
\end{equation}

\includegraphics[width=10cm, trim=3cm 6cm 5cm 4cm, clip=true]{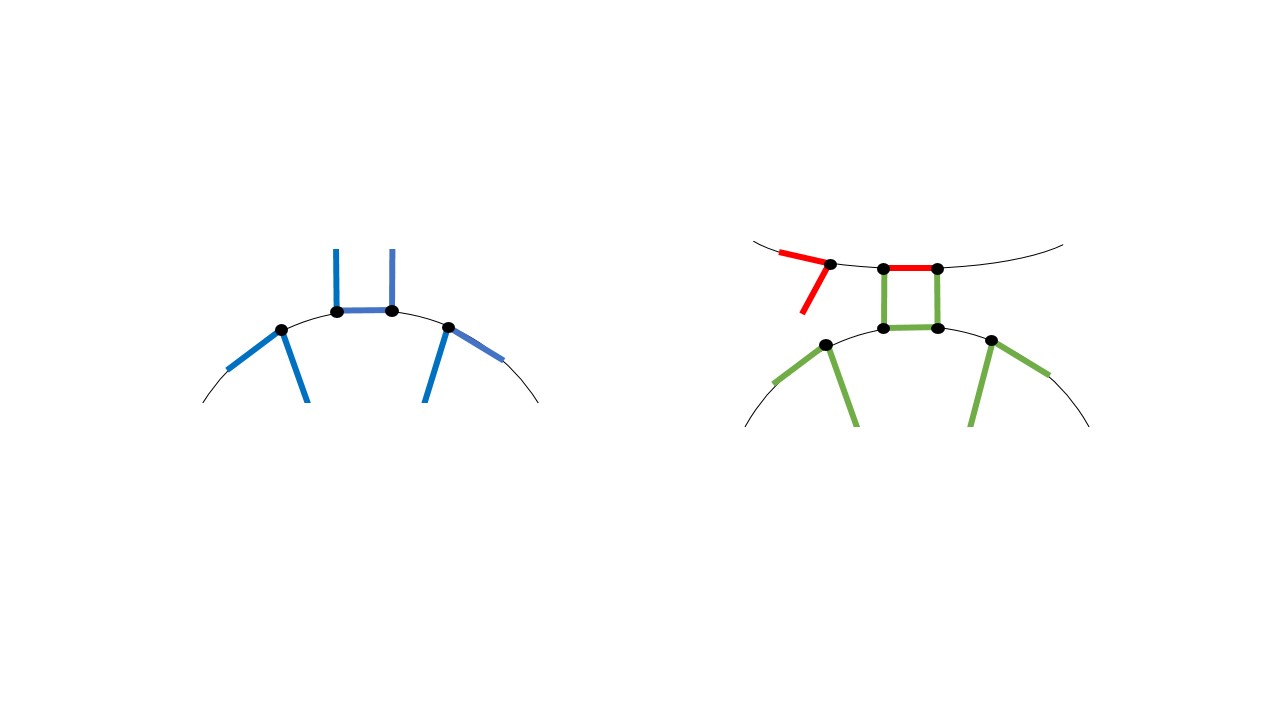}\\
Figure 3: ${\bf DC}$ (8) selects the blue edges. When the green edges are chosen, in ${\bf DC}$ (6), the red edges are chosen to close the cycle.

\includegraphics[width=9cm, trim=1cm 1cm 0cm 2cm, clip=true]{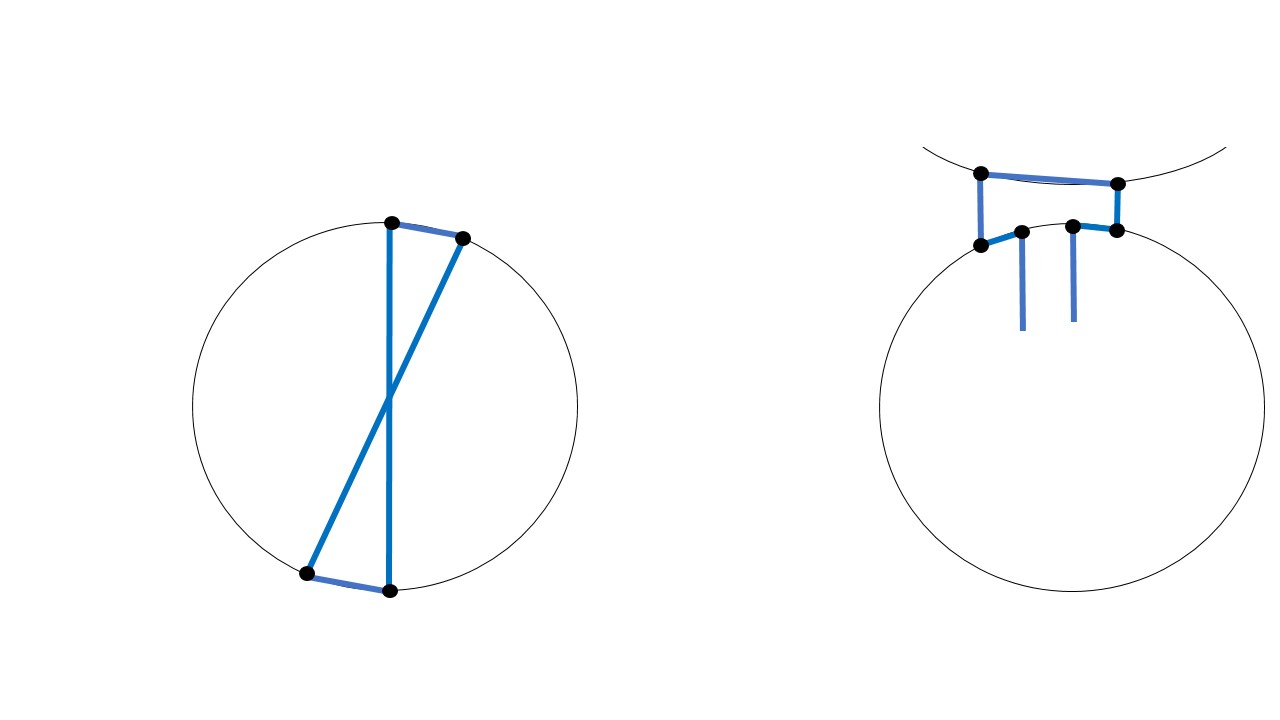}\\
Figure 4: ${\bf DC}$ (17) removes the blue cycle on the left. If the number of components of $F \cup K'$ is not changed, ${\bf DC}$ (19) may remove the blue cycle on the right.

\includegraphics[width=10cm, trim=0cm 2cm 6cm 1cm, clip=true]{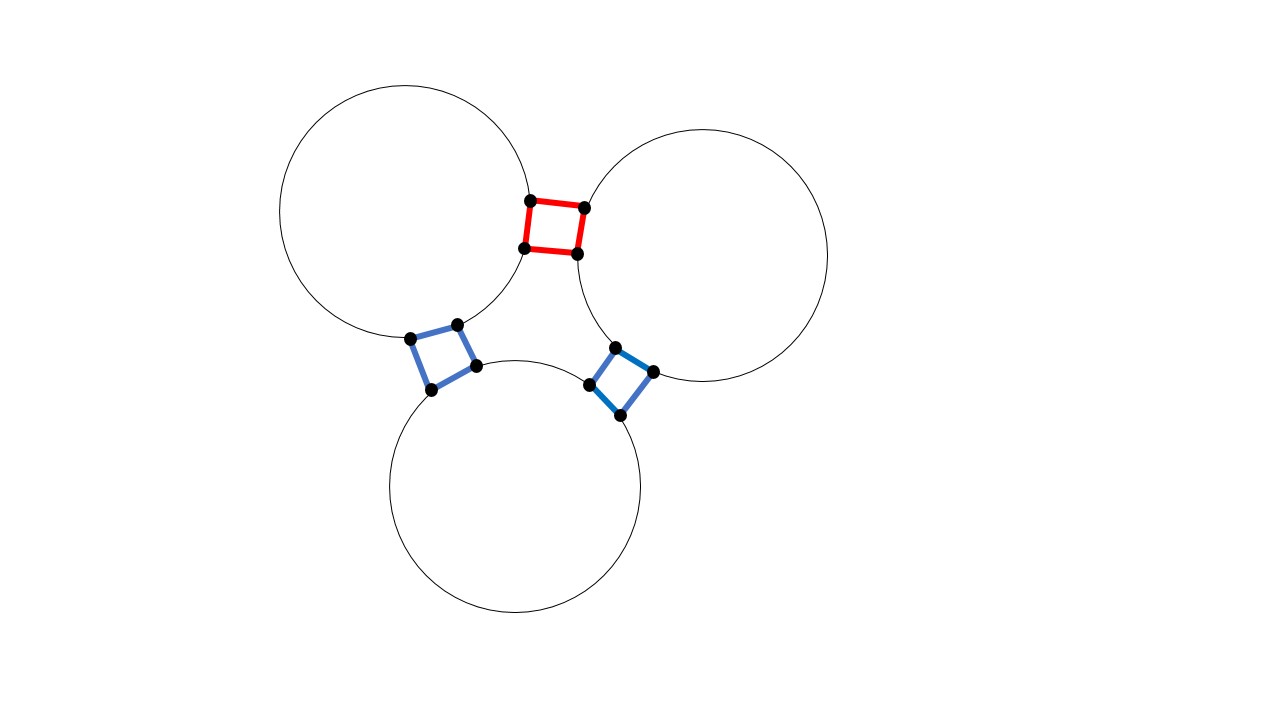}\\
Figure 5: ${\bf DC}$ (21) removes the red cycle since the number of components of $F \cup K'$ is not changed.

\begin{claim}\label{F0}
Among the cycles of $G$ with alternating edges of $F$ and edges of $M'$ or $M''$, a cycle with alternating edges of $F$ and edges of $M'$, a cycle containing a path with an edge of $M'$, an edge of $F$, an edge of $M'$ in order, and a cycle containing an edge of $M'$ and not containing a path with an edge of $M'$, an edge of $F$, an edge of $M'$ in order, are in that order, the best cycles that contribute to reducing the number of cycles in $F$ that should be included in $K$ in ${\bf HC}$ (3). Moreover, by repeating ${\bf HC}$ (3) at most $|V(G)|^2$ times, the best taken $K$ brings the number of cycles in $F$ to 1.
\end{claim}

\begin{proof}
Each cycle shown in the first half is the best to reduce the number of cycles in $F$ in ${\bf HC}$ (3) in that order because the formula (1) subtracts $F(K)$ and adds $M(K)$ to the original $F$. We check the second half. As a cycle in the best taken $K$ in ${\bf HC}$ (3) that do not contribute to reducing the number of cycles in $F$ in ${\bf HC}$ (3), a cycle containing a path with an edge of $M'$, an edge of $F$, an edge of $M'$ in order, except for cycles with alternating edges of $F$ and edges of $M'$, or a cycle containing an edge of $M'$ and not containing a path with an edge of $M'$, an edge of $F$, an edge of $M'$ in order, is considered (refer to each cycle as (a)). If the number of cycles in $F$ does not decrease in ${\bf HC}$ (3), then of the two or more cycles in the original $F$ involving (a), the cycle in which (a) takes the edge of $M''$ separates and forms a new cycle (b). Since the minimum length of a cycle that $F$ can take is 4 and $G$ is hamiltonian, (b) binds to the cycle in $F$ within $|V(G)|/4$ times of ${\bf HC}$ (3) (see Figure 6). That is, when ${\bf HC}$ (3) is repeated $n$ times ($n \geq 0$), at least $\lfloor 4n/|V(G)| \rfloor$ cycles are reduced, so the number of cycles in $F$ results in at most $|V(G)|/4 - \lfloor 4n/|V(G)| \rfloor$. If $n = |V(G)|^2$, then $|V(G)|/4 - \lfloor 4n/|V(G)| \rfloor = - 15|V(G)|/4 < 0$, so it is enough to repeat ${\bf HC}$ (3) at most $|V(G)|^2$ times.
\end{proof}

\includegraphics[width=10cm, trim=0cm 3cm 6cm 3cm, clip=true]{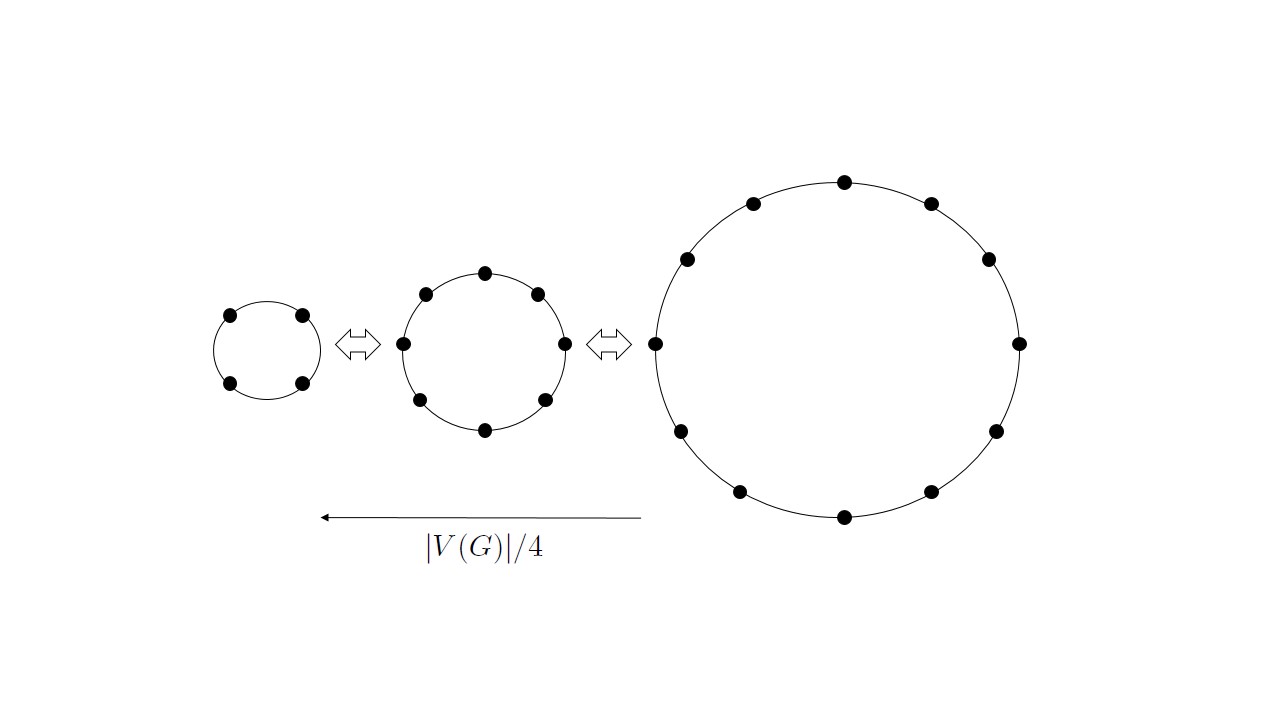}\\
Figure 6: If $G$ is hamiltonian, (b) binds to the cycle in $F$ within $|V(G)|/4$ times of ${\bf HC}$ (3). \\

Let $F_{n}$ be $F$ derived after $n$ times of ${\bf HC}$ (3), and $c_{n}$ be the number of cycles in $F_{n}$ ($n \geq 0$). By Claim \ref{F0}, $K$ is the best taken to satisfy $c_{n} = 1$ for some $n$ in $0 \leq n \leq |V(G)|^2$.
\end{proof}

\begin{thm}\label{P1}
$G$ is hamiltonian if and only if ${\bf HC}(G)$ returns a hamilton cycle. All steps of ${\bf HC}$ (including ${\bf DC}$) are computed in polynomial time. 
\end{thm}

\begin{proof}
By Lemma \ref{F1} and the definitions of ${\bf HC}$ and ${\bf DC}$, this statement obviously holds.
\end{proof}

\end{document}